\documentclass[10pt]{amsart}

\usepackage{amsthm,amsfonts,amsmath,amssymb,verbatim,cite}

\usepackage{hyperref}

\DeclareMathOperator{\height}{ht}
\DeclareMathOperator{\eup}{e}
\renewcommand{\Re}{\textup{Re}}
\newcommand{\MS}[3]{\bigg(\genfrac{}{}{0pt}{}{#1}{#2};#3\bigg)}
\newcommand{\Z}{\mathbb Z}
\newcommand{\Complex}{\mathbb C}
\newcommand{\abs}[1]{\lvert#1\rvert}
\newcommand{\ip}[2]{(#1,#2)}

\begin{document}

\begin{center}
\textbf{Ramanujan's $_1\psi_1$ summation}
\medskip

S. Ole Warnaar
\end{center}

\bigskip

\setcounter{section}{1}

\noindent \textbf{Acknowledgements}
I thank Dick Askey, Bruce Berndt, Susanna Fishel, Jeff Lagarias 
and Michael Schlosser for their helpful correspondence.

\smallskip
\noindent \textbf{Notation.}
It is impossible to give an account of the $_1\psi_1$ summation without
introducing some $q$-series notation. To keep the presentation as
simple as possible, we assume that $0<q<1$.
Suppressing $q$-dependence, we define two $q$-shifted factorials: 
$(a)_{\infty}:=\prod_{k=0}^{\infty}(1-aq^k)$ 
and $(a)_z:=(a)_{\infty}/(aq^z)_{\infty}$ for $z\in\Complex$.
Note that $1/(q)_n=0$ if $n$ is a negative integer.
For $x\in\Complex-\{0,-1,\dots\}$, the $q$-gamma function is 
defined as $\Gamma_q(x):=(q)_{x-1}/(1-q)^{x-1}$. 

\smallskip
\noindent \textbf{Ramanujan's $_1\psi_1$ summation.}
Ramanujan recorded his now famous $_1\psi_1$ summation as item~17 
of Chapter~16 in the second of his three notebooks 
\cite[p. 32]{Berndt91}, \cite{Ramanujan57}.
It was brought to the attention of the wider mathematical community in 
1940 by Hardy, who included it in his twelfth and final lecture
on Ramanujan's work \cite{Hardy40}.
Hardy remarked that the result constituted ``a remarkable formula
with many parameters''.
Instead of presenting the $_1\psi_1$ sum as given by
Ramanujan and Hardy, we will state its modern form:
\begin{equation}\label{Eq_1psi1}
\sum_{n=-\infty}^{\infty} \frac{(a)_n}{(b)_n}\, z^n=
\frac{(az)_{\infty}(q/az)_{\infty}(b/a)_{\infty}(q)_{\infty}}
{(z)_{\infty}(b/az)_{\infty}(q/a)_{\infty}(b)_{\infty}},
\qquad \abs{b/a}<\abs{z}<1,
\end{equation}
where it is understood that $a,q/b\not\in\{q,q^2,\dots\}$.
Characteristically, Ramanujan did not provide a proof of 
\eqref{Eq_1psi1}. Neither did Hardy, who however remarked that
it could be ``deduced from one which is familiar and probably goes 
back to Euler''. The result to which Hardy was referring is another famous
identity---known as the $q$-binomial theorem---corresponding to
\eqref{Eq_1psi1} with $b=q$:
$\sum_{n=0}^{\infty} z^n (a)_n/(q)_n=
(az)_{\infty}/(z)_{\infty}$ and valid for $\abs{z}<1$.
Although not actually due to Euler, the $q$-binomial theorem
is certainly classic. It seems to have appeared first and without proof
(for $a=q^{-N}$) in Rothe's 1811 book
``\textit{Systematisches Lehrbuch der Arithmetik}'', and in
the 1840s many mathematicians of note, such as Cauchy (1843), 
Eisenstein (1846), Heine (1847) and Jacobi (1847) published proofs.
The first proof of the $_1\psi_1$ sum is due to Hahn 
in 1949 \cite{Hahn49} and, as hinted by Hardy,
uses the $q$-binomial theorem.
After Hahn a large number of alternative proofs of \eqref{Eq_1psi1}
were found, including one probabilistic and three combinatorial proofs
\cite{ABBW85,Andrews69,AA78,Chan05,CCFZ11,CLM08,Chu06,CX09,Corteel03,CL3,
Ismail77,Jackson50,Kadell87,Mimachi88,Schlosser05,Yee04,VC12}.
The proof from the book, which again relies 
on the $q$-binomial theorem, was discovered by Ismail \cite{Ismail77}
and is short enough to include here. 
Assuming $\abs{z}<1$ and $\abs{b}<\min\{1,\abs{az}\}$, 
both sides are analytic functions of $b$. 
Moreover, they coincide when $b=q^{k+1}$ with $k=0,1,2,\dots$
by the $q$-binomial theorem with $a\mapsto aq^{-k}$.
Since $0$ is the accumulation point of this sequence of $b$'s
the proof is done.

Apart from the $q$-binomial theorem,
the $_1\psi_1$ sum generalises another classic identity, 
known as the Jacobi triple-product identity:
$
\sum_{n=-\infty}^{\infty} (-1)^n z^n q^{\binom{n}{2}}
=(z)_{\infty}(q/z)_{\infty}(q)_{\infty}=:\theta(z)$.
This result plays a central role in the theory of theta and
elliptic functions.

\smallskip
\noindent \textbf{The $_1\psi_1$ sum as discrete beta integral.}
As pointed out by Askey \cite{Askey78,Askey80}, the $_1\psi_1$ summation may 
be viewed as a discrete analogue of Euler's beta integral.
First define the Jackson or $q$-integral
$\int_0^{c\cdot\infty} f(t) \textup{d}_q t:=
(1-q)\sum_{n=-\infty}^{\infty} f(cq^n) cq^n$.
Replacing $(a,b,z)\mapsto (-c,-cq^{\alpha+\beta},q^{\alpha})$ 
in \eqref{Eq_1psi1} then gives
\begin{equation}\label{Eq_qint}
\int_0^{c\cdot\infty} \frac{t^{\alpha-1}}{(-t)_{\alpha+\beta}}\,
\textup{d}_q t = c^{\alpha} \,
\frac{\theta(-cq^{\alpha})}{\theta(-c)}\,
\frac{\Gamma_q(\alpha)\Gamma_q(\beta)}{\Gamma_q(\alpha+\beta)},
\end{equation}
where $\Re(\alpha),\Re(\beta)>0$.
For real, positive $c$ the limit $q\to 1$ can be taken,
resulting in the beta integral modulo the substitution $t\mapsto t/(1-t)$. 
Askey further noted in \cite{Askey78} that the specialisation 
$(\alpha,\beta)\mapsto (x,1-x)$ in \eqref{Eq_qint} (so that $0<\Re(x)<1$)
may be viewed as a $q$-analogue of Euler's reflection formula.

\smallskip
\noindent \textbf{Simple applications of the $_1\psi_1$ sum.}
There are numerous easy applications of the $_1\psi_1$ sum.
For example, Jacobi's well-known four- and six-square theorems 
as well as a number of similar results 
readily follow from \eqref{Eq_1psi1}, see e.g., 
\cite{Adiga92,BA88,Chan04,Cooper01,CL02,Fine88}. 
To give a flavour of how the $_1\psi_1$ implies these types of results
we shall sketch a proof of the four-square theorem.
Let $r_s(n)$ be the number of representations of $n$ as the sum of $s$ 
squares. The generating function $R_s(q):=\sum_{n\geq 0} r_s(n) (-q)^n$ 
is given by $\big(\sum_{m=-\infty}^{\infty} (-1)^m q^{m^2}\big)^s$. 
By the triple-product identity this is also 
$\big((q)_{\infty}/(-q)_{\infty}\big)^s$.
Any identity that allows the extraction of
the coefficient of $(-q)^n$ results in an explicit formula
for $r_s(n)$. Back to \eqref{Eq_1psi1}, replace 
$(b,z)\mapsto (aq,b)$ and multiply both sides by 
$(1-b)/(1-ab)$. By the geometric series this yields
\begin{equation}\label{Eq_Kronecker}
1+\frac{(1-a)(1-b)}{1-ab}\sum_{k,n=1}^{\infty} q^{kn}(a^k b^n-a^{-k}b^{-n})
=\frac{(abq)_{\infty}(q/ab)_{\infty}(q)_{\infty}^2}
{(aq)_{\infty}(q/a)_{\infty}(bq)_{\infty}(q/b)_{\infty}},
\end{equation}
which may also be found in Kronecker's 1881 paper 
``\textit{Zur Theorie der elliptischen Functionen}''.
For $a,b\to -1$ the right side gives $R_4(q)$ whereas the left side becomes
\[
1-8\sum_{m=1}^{\infty}q^m 
\sum_{\substack{n,k=1 \\ nk=m}}^{\infty} n (-1)^{n+k}=
1+8\sum_{m=1}^{\infty}(-q)^m 
\sum_{\substack{d\geq 1 \\ 4 \nmid d\mid m}} d.
\]
Hence $r_4(n)=8\sum_{d\geq 1;~4 \nmid d\mid n} d$. This result
of Jacobi implies Lagrange's theorem that every positive integer is
a sum of four squares.
By taking $a,b^2\to -1$ in \eqref{Eq_Kronecker} the reader will have little 
trouble showing that $r_2(n)=4(d_1(n)-d_3(n))$,
with $d_k(n)$ the number of divisors of $n$ of the form $4m+k$.
This is a result of Gauss and Lagrange which implies Fermat's
two-square theorem.

Other simple but important applications of the $_1\psi_1$ sum 
concern orthogonal polynomials. In \cite{AW85} it was employed 
by Askey and Wilson to compute a special 
case---corresponding to the continuous $q$-Jacobi polynomials---of the 
Askey--Wilson integral, and in \cite{Askey83} Askey gave an elementary
proof of the full Askey--Wilson integral using the $_1\psi_1$ sum. 
The sum also implies the norm evaluation of the weight 
functions of the %$q$-ultraspherical polynomials \cite{AI83} and 
$q$-Laguerre polynomials \cite{Moak81}. 
These are a family of orthogonal polynomials
with discrete measure $\mu$ on $[0,c\cdot\infty)$ given by
$\textup{d}_q \mu(t)=t^{\alpha}/(-t)_{\infty}\, \textup{d}_q t$.
The normalisation $\int \textup{d}_q \mu(t)$ thus follows from the
$q$-beta integral \eqref{Eq_qint} in the limit of large $\beta$.

\smallskip
\noindent \textbf{Generalisations in one dimension.}
There exist several generalisations of Ramanujan's sum containing
one additional parameter. In his work on partial theta functions 
Andrews \cite{Andrews81} obtained a generalisation in which each product of
four infinite products on the right-hand side is replaced by six such products. 
Another example is the curious identity of Guo and 
Schlosser, which is no longer hypergeometric in nature \cite{GS07}:
\[
\sum_{k=-\infty}^\infty 
%\frac{(a)_k}{(b)_k}\,
%\frac{(1-ac_kq^k)}{(1-azq^k)}\,
%\frac{(c_kq)_{\infty}(b/ac_k)_{\infty}}
%{(ac_k\big)_{\infty}(q/ac_k\big)_{\infty}}
\frac{(a)_k (1-ac_kq^k)(c_kq)_{\infty}(b/ac_k)_{\infty}}
{(b)_k (1-azq^k)(ac_k\big)_{\infty}(q/ac_k\big)_{\infty}}\,
c_k^k
=\frac{1}{(1-z)}\,
\frac{(q)_\infty(b/a)_\infty}{(q/a)_\infty(b)_\infty},
\]
where $c_k:=z(1-aczq^k)/(1-azq^k)$ and 
$|b/ac|<\abs{z}<1$. For $c=1$ this is \eqref{Eq_1psi1}.

As discovered by Schlosser \cite{Schlosser06},
a quite different extension of the $_1\psi_1$ sum arises
by considering non-commutative variables.
Let $R$ be a unital Banach algebra with identity $1$,
central elements $b$ and $q$, and norm $\|\cdot\|$. Write
$a^{-1}$ for the inverse of an invertible element $a\in R$.
Let $\prod_{i=m}^n a_i$ stand for $1$ if $n=m-1$, $a_m\cdots a_n$
if $n\geq m$ and $a_{m-1}^{-1}\cdots a_{n+1}^{-1}$ if $n<m-1$, and
define
\[
\MS{a_1,\dots,a_r}{b_1,\dots,b_r}{z}_{\! k}^{\!\pm}:=
\prod_i \bigg[z \prod_{s=1}^r 
(1-a_s q^{i-1})(1-b_s q^{i-1})^{-1}\bigg],
\]
where $k\in\Z\cup\{\infty\}$, $a_1,\dots,a_r,b_1,\dots,b_r\in R$,
$\prod_i=\prod_{i=1}^k$ in the $+$ case and
$\prod_i=\prod_{i=k}^1$ in the $-$ case. 
Subject to $\max\{\|q\|,\|z\|,\|ba^{-1}z^{-1}\|\}<1$,
the following non-commutative $_1\psi_1$ sum holds:
\[
\sum_{k=-\infty}^{\!\infty}
\MS{a}{b}{z}^{\! +}_{\! k}=
\MS{za}{z}{1}^{\! -}_{\!\infty}
\MS{qa^{-1}z^{-1}}{qza^{-1}z^{-1}}{1}^{\! +}_{\!\infty}
\MS{bza^{-1}z^{-1},q}{ba^{-1}z^{-1},b}{1}^{\! -}_{\!\infty}.
\]

\smallskip
\noindent \textbf{Higher-dimensional generalisations.}
Various authors have generalised \eqref{Eq_1psi1} to multiple 
$_1\psi_1$ sums. Below we state a generalisation 
due to Gustafson and Milne \cite{Gustafson87,Milne86} which is
labelled by the $\mathrm{A}$-type root system.
Similar such $_1\psi_1$ sums are given in
\cite{Aomoto95,Aomoto98,GK96,MS02,Schlosser00}.
More involved multiple $_1\psi_1$ sums with a Schur or Macdonald polynomial 
argument can be found in \cite{BF99,Kaneko96,Milne92,Warnaar10}.
For $r=(r_1,\dots,r_n)\in\Z^n$ denote $\abs{r}:=r_1+\cdots+r_n$. Then
\[
\sum_{r\in\Z^n} z^{\abs{r}} \prod_{1\leq i<j\leq n} 
\frac{x_i q^{r_i}-x_j q^{r_j}}{x_i-x_j}\,
\prod_{i,j=1}^n\frac{(a_j x_{ij})_{r_i}}{(b_j x_{ij})_{r_i}}
=\frac{(az)_{\infty}(q/az)_{\infty}}{(z)_{\infty}(b/az)_{\infty}} 
\prod_{i,j=1}^n \frac{(b_jx_{ij}/a_i)_{\infty} (qx_{ij})_{\infty}}
{(qx_{ij}/a_i)_{\infty}(b_j x_{ij})_{\infty}},
\]
where $a:=a_1\cdots a_n$, $b:=q^{1-n}b_1\cdots b_n$, $x_{ij}:=x_i/x_j$
and $\abs{b/a}<\abs{z}<1$. Milne first proved this for $b_1=\dots=b_n$ 
\cite{Milne86} and shortly thereafter Gustafson established the full
result \cite{Gustafson87}.
We have already seen that the $_1\psi_1$ sum implies the 
Jacobi triple-product identity. The latter is the $\mathrm{A}_1^{(1)}$
case of Macdonald's generalised Weyl denominator identities
for affine root systems \cite{Macdonald72}. 
To obtain further Macdonald identities 
from the Gustafson--Milne sum one replaces $z\to z/a$ before letting
$a_1,\dots,a_n\to\infty$ and $b_1, \dots,b_n\to 0$. 
Extracting the coefficient of $z^0$ (on the right this
requires the triple-product identity) results in the 
Macdonald identity for $\mathrm{A}_{n-1}^{(1)}$.

Higher-dimensional generalisations of a special case of the $_1\psi_1$ sum 
can be given for all affine root systems. A full description is beyond this 
note, and we will only sketch the simplest case.
The reader is referred to \cite{FGT08,Macdonald72,Macdonald72b,Macdonald03}
for the full details.
In \cite{Macdonald72b} Macdonald gave the following multivariable 
extension of the product formula for the Poincar\'e polynomial of a 
Coxeter group
\begin{equation}\label{Eq_Macdonald}
W(\boldsymbol{t}):=
\sum_{w\in W} \prod_{\alpha\in R^{+}}
\frac{1-t_{\alpha} \eup^{w(\alpha)}}
{1-\eup^{w(\alpha)}}
=\prod_{\alpha\in R^{+}} \frac{1-t_{\alpha}\boldsymbol{t}^{\height(\alpha)}}
{1-\boldsymbol{t}^{\height(\alpha)}}.
\end{equation}
Here $R$ is a reduced, irreducible finite root system in a Euclidean 
space $V$, $R^{+}$ the set of positive roots, 
$W$ the Weyl group and $t_{\alpha}$ 
for $\alpha\in R^{+}$ a set of formal variables constant along Weyl orbits.
The symbol $\boldsymbol{t}^{\height(\alpha)}$ stands for
$\prod_{\beta\in R^{+}} t_{\beta}^{\ip{\beta}{\alpha}/\|\alpha\|^2}$
with $\ip{\cdot}{\cdot}$ the $W$-invariant positive 
definite bilinear form on $V$. If all $t_{\alpha}$ are set to $t$ then 
$\boldsymbol{t}^{\height(\alpha)}=t^{\height(\alpha)}$ with $\height(\alpha)$
the usual height function on $R$, in which case $W(\boldsymbol{t})$ reduces to
the classical Poincar\'e polynomial $W(t)$.
Now let $S$ be a reduced, irreducible affine root
system of type $S=S(R)$ \cite{Macdonald72}. In analogy with the finite case,
assume that $t_a$ for $a\in S$ is constant along orbits of the 
affine Weyl group $W$ of $S$. Then Macdonald generalised \eqref{Eq_Macdonald}
to \cite{Macdonald03}
\begin{equation}\label{Eq_Macdonald_b}
\sum_{w\in W} \prod_{a\in S^{+}}
\frac{1-t_a \eup^{w(a)}}
{1-\eup^{w(a)}}
=\prod_{\alpha\in R^{+}}
\frac{(t_{\alpha}\boldsymbol{t}^{\height(\alpha)})_{\infty}
(\boldsymbol{t}^{\height(\alpha)}q^{\chi(\alpha\in B)}/t_{\alpha})_{\infty}}
{(\boldsymbol{t}^{\height(\alpha)})_{\infty}^2},
\end{equation}
where $B$ is a base for $R$.
The parameter $q$ on the right is fixed by $q=\prod_{a\in B(S)} \exp(n_a a)$,
where $B(S)$ is a basis for $S$ and the $n_a$ are the labels of the 
extended Dynkin diagrams given in \cite{Macdonald72}.
If $R$ is simply-laced then $t_a=t$.
In the case of $S(R)=\textrm{A}_1^{(1)}$,
$q=\exp(a_0+a_1)$ so that after replacing $\exp(a_1)$ by $x$ we obtain
the $_1\psi_1$ sum \eqref{Eq_1psi1} with $(a,b,z)\to (x/t,tx,t)$.
This is not the end of the story concerning root systems and the
$_1\psi_1$ sum. Identity \eqref{Eq_Macdonald_b} can be rewritten 
as \cite{Macdonald03}
\begin{equation}\label{Eq_Macdonald_c}
\sum_{\gamma\in Q^{\vee}} \prod_{\alpha\in R}
\frac{(q\eup^{\alpha})_{\ip{\alpha}{\gamma}}}
{(t_{\alpha} q\eup^{\alpha})_{\ip{\alpha}{\gamma}}}
=\prod_{\alpha\in R^{+}}
\frac{(t_{\alpha}\boldsymbol{t}^{\height(\alpha)}q)_{\infty}
(\boldsymbol{t}^{\height(\alpha)}q^{\chi(\alpha\in B)}/t_{\alpha})_{\infty}}
{(\boldsymbol{t}^{\height(\alpha)}q)_{\infty}
(\boldsymbol{t}^{\height(\alpha)})_{\infty}}\,
\frac{(q\eup^{\alpha})_{\infty}(q\eup^{-\alpha})_{\infty}}
{(t_{\alpha} q\eup^{\alpha})_{\infty}
(t_{\alpha} q\eup^{-\alpha})_{\infty}},
\end{equation}
where $Q^{\vee}$ is the coroot lattice.
Interestingly, for $t_{\alpha}=t$ this was also found by 
Fishel, Grojnowski and Teleman \cite{FGT08} by computing
the generating function of the $q$-weighted Euler 
characteristics of certain Dolbeault cohomologies. 
For $R=\mathrm{A}_{n-1}$, $Q^{\vee}=Q=
\sum_{i=1}^n r_i\epsilon_i$ with $\abs{r}=0$, 
$R=\{\epsilon_i-\epsilon_j:~1\leq i\neq j\leq n\}$ and 
$t^{\height(\epsilon_i-\epsilon_j)}=t^{j-i}$.
By fairly elementary manipulations the identity \eqref{Eq_Macdonald_c}
may then be transformed into the multiple $_1\psi_1$ sum
\begin{multline*}
\sum_{r\in\Z^n}  \: z^{\abs{r}}\frac{(a)_{\abs{r}}}{(b)_{\abs{r}}}
\prod_{1\leq i<j\leq n} \frac{x_iq^{r_i}-x_jq^{r_j}}{x_i-x_j}\,
\frac{(t^{-1}x_{ij})_{r_i-r_j}}{(tqx_{ij})_{r_i-r_j}}
\, t^{r_i-r_j}q^{-r_j} \\
=\frac{(az)_{\infty}(q/az)_{\infty}(b/a)_{\infty}(tq)_{\infty}}
{(z)_{\infty}(b/az)_{\infty}(q/a)_{\infty}(b)_{\infty}}\,
\prod_{i=1}^{n-1} \frac{(t^{i+1}q)_{\infty}}{(t^i)_{\infty}}
\prod_{i,j=1}^n \frac{(qx_{ij})_{\infty}}{(tqx_{ij})_{\infty}},
\end{multline*}
for $\abs{b/a}<\abs{z}<1$ and $\abs{t}<1$.
This is the only result in this survey that is new.

We finally remark that all higher-dimensional $_1\psi_1$ sums
admit representations as discrete Selberg-type integrals.
The most important such integrals are due to Aomoto
\cite{Aomoto95,Aomoto98} and Ito \cite{Ito97}, and are closely related to
\eqref{Eq_Macdonald_b}. Further examples may be found in
\cite{Kaneko96b,Warnaar05}.

\small

\bibliographystyle{amsplain}

\medskip

\tiny {

\noindent
School of Mathematics and Physics, The University of Queensland, 
Brisbane, QLD 4072, Australia \\
Work supported by the Australian Research Council
}

\end{document}